\documentclass[11pt]{article}

\usepackage{amsmath, amssymb, amsthm, graphicx}

\usepackage{amsfonts}

\usepackage{booktabs}

\usepackage{float}

\title{Representability of Lyndon-Maddux relation algebras}
\author{Jeremy F.~Alm}
\date{March 2017}

\newtheorem{thm}{Theorem}
\newtheorem{defn}[thm]{Definition}

\def\rp{{\; ; \;}}
\def\id{{1\kern-.08em\raise1.3ex\hbox{\rm,}\kern.08em}}
\def\di{{0\kern-.04em\raise1.3ex\hbox{,}\kern.04em}}

\def\L{\mathfrak{L}}
\def\RA{{\sf RA}}
\def\RRA{{\sf RRA}}
\def\wRRA{{\sf wRRA}}
\def\Pr{{\sf Pr}}

\begin{document}

\maketitle

\begin{abstract}
    In Alm-Hirsch-Maddux (2016), relation algebras $\mathfrak{L}(q,n)$ were defined that generalize Roger Lyndon's relation algebras from projective lines, so that $\mathfrak{L}(q,0)$ is a Lyndon algebra.  In that paper, it was shown that if $q>2304n^2+1$, $\mathfrak{L}(q,n)$ is representable, and if $q<2n$, $\mathfrak{L}(q,n)$ is not representable.  In the present paper, we reduced this gap by proving that if $q\geq n(\log n)^{1+\varepsilon}$,  $\mathfrak{L}(q,n)$ is representable.
\end{abstract}

\section{Introduction}

Let \RA\ denote the class of relation algebras,  let
\RRA\ denote the class of  representable relation algebras, and let
\wRRA\ denote the class of weakly representable relation algebras.  Then we have $\RRA \subsetneq \wRRA \subsetneq\RA$, and \RRA\ (resp., \wRRA) is not finitely axiomatizable over \wRRA\ (resp., \RA).  J\'onsson \cite{Jon91} showed that  every equational basis for \RRA\ contains equations with arbitrarily many variables. Since all three of \RA, \wRRA, and \RRA\ are varieties, we may infer from J\'onsson's result either that   every equational basis for \wRRA\ contains equations with arbitrarily many variables or every equational basis defining \RRA\ over \wRRA\ contains equations with arbitrarily many variables.  In \cite{AHM}, it is shown that the latter holds. (Whether the former holds is still open.)

The argument in \cite{AHM}  uses a construction of
arbitrarily large finite weakly representable but not representable
relation algebras whose ``small'' subalgebras are representable.  Define a class of algebras $\L(q,n)$ as follows.
\begin{defn}\label{def:madd}
Let $q$ be a positive integer, and $n$ a nonnegative integer. 
Then $\L(q,n)$ is the symmetric integral relation algebra with   atoms $\id, a_0,\cdots,a_q,  t_1,\cdots,t_n$. 
Write $A=a_0+\cdots+a_p$ and $T=t_1+\cdots+t_n$, and define $\rp$ on atoms
as follows: if $0\leq i,j\leq q$, $i\neq j$, $1\leq k,l\leq n$, and
$k\neq l$, then
\begin{align*}
a_i\rp a_i&=\id+a_i,\\
a_i\rp a_j&=A\cdot\overline{a_i+a_j},\\
a_i\rp t_k&=T,\\
t_k\rp t_k&=\id+A,\\
t_k\rp t_l&=A.
\end{align*}
\end{defn}

Definition \ref{def:madd} from \cite{AHM} is due to Maddux.  When $n=0$, $\L(q,0)$ is the Lyndon algebra from a  projective line. $\L(q,0)$ is representable over $q^2$ points exactly when there exists a projective plane of order $q$.  For $n=1$, the idea is to take a Lyndon algebra $\L(q,0)$, add an additional atom $t$ and mandatory cycles $att$, where $a$ is any atom besides $t$, to get $\L(q,1)$. For $n>1$, split $t$ (in the sense of \cite{AMN}) into $n$ smaller atoms. Even if $\L(q,1)$ is representable, splitting $t$  into too many pieces can destroy representability (although in this particular case splitting does preserve weak representability). In fact, we have the following theorem from \cite{AHM}.

\begin{thm} \label{thm:bound} Let $q$ be a prime power, so that $\L(q,0)$ is representable. Then 
\begin{enumerate}
\item[(i.)] if $2n>q$, then $\L(q,n)\notin\RRA$;
\item[(ii.)] If $q>2304n^2+1$, then $\L(q,n)\in\RRA$.
\end{enumerate}
\end{thm}

Our concern in the present paper is this: How large does $q$ need to be, relative to $n$, to guarantee that if $q$ is a prime power, $\L(q,n)\in\RRA$?  The sole result of the present paper is an improvement of the bound in Theorem \ref{thm:bound} (ii.), as follows.
\begin{thm}
Let $\varepsilon>0$, and let $q$ be a  prime power with $q\geq n(\log n)^{1+\varepsilon}$. Then for sufficiently large $n$,  $\L(q,n)\in\RRA$.  In particular, for $n\geq 18$, if $q\geq n(\log n)^3$, $\L(q,n)\in\RRA$.
\end{thm}

\section{The union bound and the Lov\'asz local lemma}

In this section we will establish the necessary background in probabilistic combinatorics.  We will do this somewhat informally, since there are plenty of rigorous references on the subject, like \cite{AlonSpencer, Jukna}.

Consider a random combinatorial structure $R$, along with a probability space corresponding to all possible particular instances of $R$. (For example, for the standard random graph model $G_{n,0.5}$, where edges are either included or excluded independently with probability 0.5, the probability space would be the space of all graphs with $n$ vertices.) Suppose that we want to show that $R$ can have a property $\mathcal{P}$. Suppose further that is a list of events $\{A_i\}_{i\in I}$, the ``bad'' events, that prevent $\mathcal{P}$ from being attained. (For example, if $\mathcal{P}$ is the property of being bipartite, then the events $A_i$ are the instances of odd cycles.)  We would like to estimate the probability $\Pr\left[\bigcup A_i\right]$ that any of the bad events occur.  In particular, in many cases one may show $\Pr\left[\bigcup A_i\right] \to 0$ as the size of a random combinatorial structure grows, but all we really need for a nonconstructive existence proof is to show $\Pr\left[\bigcup A_i\right] < 1$.
The following two tools are widely used to this end.

\textbf{Union Bound}

For events $A_i$,
\[
    \Pr\left[\bigcup A_i\right]\leq \sum \Pr[A_i].
\]

Informally, we may think of the union bound as a simpleton's inclusion-exclusion principle, whereby we throw away all the terms involving intersections. 

\textbf{Lov\'asz local lemma} 

If the $A_i$ each occur with probability at most $p<1$, and with each event independent of all other $A_i$'s except at most $d$ of them, then if $epd\leq 1$, then $\Pr\left[\bigcup A_i\right] < 1$. (Here, $e$ is the base of the natural logarithm.)



\section{Proof of main result}

Consider $\mathbb{F}_q \times \mathbb{F}_q$, where  $\mathbb{F}_q$ is the field with $q$ elements. For convenience, let $\mathbb{F}_q = \{ 0, 1, \dots, q-1 \}$.   For each $i \in \mathbb{F}_q$, define
\[
    R_i = \{ \big( (a_1,b_1),(a_2,b_2) \big) \in (\mathbb{F}_q \times \mathbb{F}_q)^2 : (b_2-b_1)\cdot(a_2-a_1)^{-1} = i\}
\]
and
\[
    R_q = \{ \big( (a_1,b_1),(a_2,b_2) \big) \in (\mathbb{F}_q \times \mathbb{F}_q)^2 : a_2=a_1 \}
\]
One can think of $R_i$ as the union of lines with slope $i$, and of $R_q$ as the union of vertical lines, in $\mathbb{F}_q \times \mathbb{F}_q$.  As is well-known, the map $a_i \mapsto R_i$ is a representation of $\L(q,0)$. (This was originally shown in \cite{Lyndon}.)  Write $D=\mathbb{F}_q \times \mathbb{F}_q$.  Let $D'$ be a disjoint copy of $D$, and let $R_i' \subset D' \times D'$ be the corresponding copy of $R_i$.  Then the map that sends $a_i \mapsto R_i \cup R_i'$ and $t$ to $D \times D' \cup D' \times D$ is a representation of $\L(q,1)$. (See \cite{AHM}.)

Now we define a random structure $R(q,n)$. Take the (image of the) representation of $\L(q,1)$, and view it as a complete graph with vertex set $D\cup D'$. Consider all possible ways of coloring the $t$-edges in colors $t_1,\ldots ,t_n$. We show that the resulting probability space contains a representation of $\L(q,n)$.

We need to check that all of the following conditions obtain:
\begin{equation} \label{att} 
    \text{for every } k\leq q+1,
 \text{ for every } i,j\leq n,\ a_k\leq t_i\rp t_j.
\end{equation}

\begin{equation}  \label{tatta} 
\text{for every } i\leq n, j\leq n, \text{ and } k\leq q+1,\
t_i\leq a_k\rp t_j\cdot t_j\rp a_k.
\end{equation}

If there is an edge $uv$, labelled $a_k$, say, for which there is no $z$ such that $uz$ is labelled $t_i$ and $vz$ is labelled $t_j$, we will say that that edge fails. (In general, conditions \eqref{att} and \eqref{tatta} give each edge some ``needs'', and we will say that an edge fails if it has any unmet needs.)

Fix an edge labelled $a_k$.  Each vertex has $q^2$ edges colored in $t$-colors, so the probability that \eqref{att} fails on a given $a_k$-edge for colors $(t_i,t_j)$ is $(1-\frac{1}{n^2})^{q^2}$.  Hence the failure for any such pair $(t_i,t_j)$ is bounded above (via the union bound) by 
\[
    \sum_{t_i,t_j}\left(1-\frac{1}{n^2}\right)^{q^2}=n^2\left(1-\frac{1}{n^2}\right)^{q^2}.
\]

Similarly, fix an edge colored $t_i$. Each vertex has degree $q-1$ in each $a_k$ color, and degree $q^2$ in the various $t_j$ colors. Let the edge colored $a_k$ be $(u,v)$. The probability that there is no $z$ with $uz$ colored $a_k$ and $zv$ colored $t_j$ for fixed $k,j$ is $(1-\frac{1}{n})^{q-1}$. Hence the failure for all $k,j$ is bounded by $n(q+1)(1-\frac{1}{n})^{q-1}$. Since we must also consider $z'$ such that $uz'$ is colored $t_j$ and $z'v$ is colored $a_k$, we double to get $2n(q+1)(1-\frac{1}{n})^{q-1}$.  Now, we may bound the probability that there exists an edge that fails. Since there are $2\binom{q^2}{2}$ $a$-labelled edges and $q^4$ $t$-labelled edges, the union bound gives
\begin{equation} \label{eqn:binom}
    2\binom{q^2}{2}\cdot n^2\left(1-\frac{1}{n^2}\right)^{q^2} + q^4\cdot 2n(q+1)\left(1-\frac{1}{n}\right)^{q-1}
\end{equation}

\begin{figure}[hbt]
    \centering
    \includegraphics[width=4in]{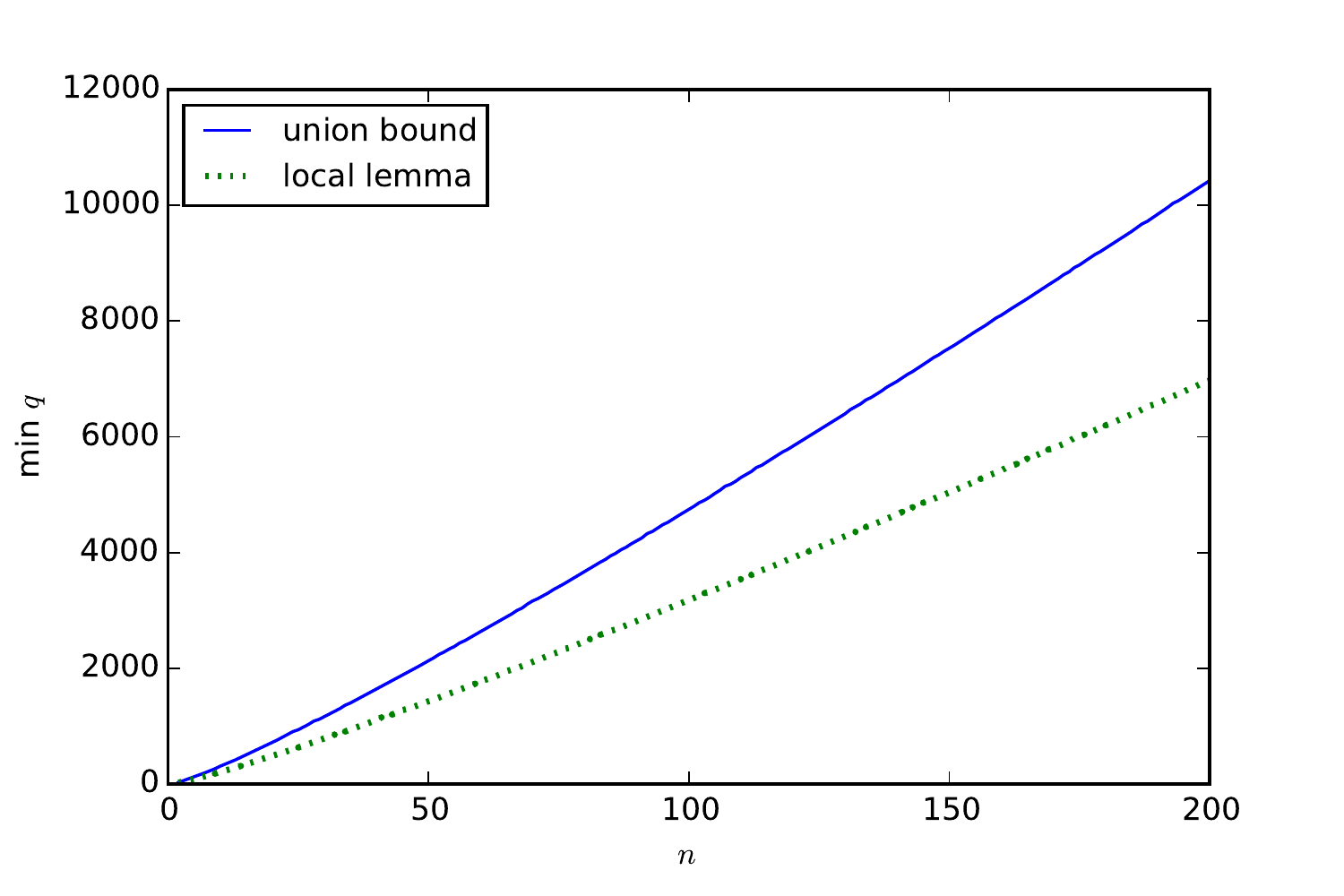}
    \caption{Smallest prime power $q$ as a function of $n$}
    \label{fig:ubvslll}
\end{figure}

We want to make \eqref{eqn:binom} $<1$, so  there exists positive probability that no edge fails.  So far, so good, but we can do better using the local lemma.  (In particular, see Table \ref{tab:qvsn} and Figure \ref{fig:ubvslll}.)

In order to make use of the local lemma, we need to count dependencies.  For  any edge $uv$, the failure of that edge is dependent on all the $a$-labelled edges incident to $u$ and $v$, as well as all of the $t$-labelled edges.  So we bound this by $4q^2$, twice the total number of vertices.

We need $edp\leq 1$. Consequently, it must be the case that
\[
    d\leq 4q^2, \text{ and } p\leq 2n\left(q+1\right)\left(1-\frac{1}{n}\right)^{q-1}. 
\]

So we want
\[
    e\cdot 4q^2\cdot 2n(q+1)\left(1-\frac{1}{n}\right)^{q-1}\leq 1.
\]

Taking logs of both sides and manipulating, we get
\begin{equation} \label{eqn:lll}
     1+\log (8)+\log(n) +2\log(q)+\log(q+1)\leq (q-1)\log\left(\frac{n}{n-1}\right).
\end{equation}

Since $\log\left(\frac{n}{n-1}\right)\sim \frac{1}{n}$, we need to take $q$ slightly larger than $n\log n$, so that the RHS will eventually beat the $\log(n)$ on the LHS.

\begin{figure}[H]
    \centering
    \includegraphics[width=4in]{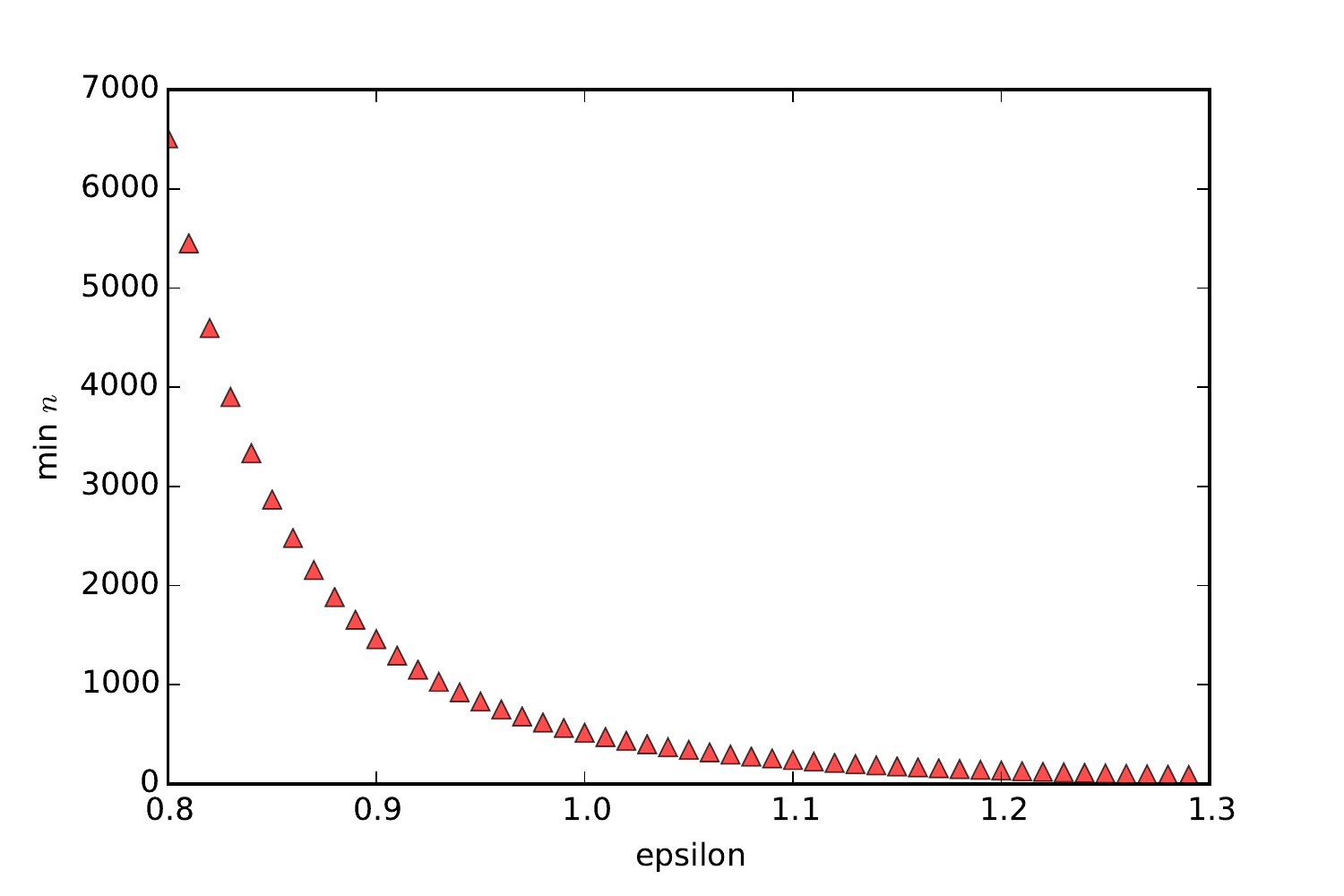}
    \caption{Smallest $n$ for which $n(\log n)^{1+\varepsilon}$ satisfies \eqref{eqn:lll}}
    \label{fig:nvsepsilon}
\end{figure}

Taking $q\geq n(\log n)^{1+\varepsilon}$, $\varepsilon>0$, will do for sufficiently large $n$ (depending on $\varepsilon$).  See Figure \ref{fig:nvsepsilon} for the dependence of $n$ upon $\varepsilon$.

\begin{table}[H]
    \centering
   \begin{tabular}{ccc}
$n$ &  union bound  & local lemma \\
\toprule
2 & 27 & 23\\
3 & 59 & 41\\
4 & 89 & 61\\
5 & 121 & 83\\
6 & 157 & 107\\
7 & 191 & 131\\
8 & 227 & 157\\
9 & 263 & 179\\
10 & 307 & 211\\
11 & 343 & 233\\
12 & 379 & 257\\
13 & 419 & 289\\
14 & 461 & 311\\
15 & 503 & 343\\
16 & 547 & 367\\
17 & 587 & 397\\
18 & 631 & 431\\
19 & 673 & 457\\
20 & 719 & 487
\end{tabular}
    \caption{Smallest prime powers for small $n$ for which \eqref{eqn:binom} $<1$, resp., \eqref{eqn:lll} holds}
    \label{tab:qvsn}
\end{table}

\bibliographystyle{plain}
\bibliography{refs}

\end{document}